\documentclass[12pt]{article}
\usepackage[utf8]{inputenc}
\usepackage[french,english]{babel}
\usepackage[T1]{fontenc}
\usepackage{amsmath}
\usepackage{amsfonts}
\usepackage{amssymb}
\usepackage{graphicx}
\usepackage[dvipsnames]{xcolor}
\usepackage{tcolorbox}
\usepackage{hyperref}
\usepackage[nottoc]{tocbibind}
\usepackage[ right=2cm, left=2cm, top=2cm, bottom=2cm]{geometry}
\usepackage{cite}

\newtheorem{Theorem}{Theorem}[section]
\newtheorem{Lemma}{Lemma}[section]
\newtheorem{Proposition}{Proposition}[section]
\newtheorem{Definition}{Definition}[section]
\newtheorem{Remark}{Remark}[section]

\newgeometry{includefoot,left=2cm,right=2cm,bottom=1cm,top=2cm}
\title{ \LARGE \bf%
{Weak input to state estimates for 2D damped wave equations with localized and non-linear damping}
\thanks{This research was partially supported by the iCODE Institute, research project of the IDEX Paris-Saclay, and by the Hadamard Mathematics LabEx (LMH) through the grant number ANR-11-LABX-0056-LMH in the ``Programme des Investissements d'Avenir''.}
}
\author{Meryem Kafnemer$^{1}$, Benmiloud Mebkhout$^{2}$ and Yacine Chitour$^{3}$}
\footnotetext[1]{UMA, ENSTA,  Palaiseau and  L2S, Universit\'e Paris Saclay, France.}
\footnotetext[2]{D\'epartement de Math\'ematiques, Université Abou Bekr Belkaid, Tlemcen, Algeria.}
\footnotetext[3]{L2S, Universit\'e Paris Saclay, France.}
\begin{document}
\maketitle
\begin{abstract}
	\ \\
In this paper, we study input-to-state (ISS) issues for damped wave equations with Dirichlet boundary conditions on a bounded domain of dimension two. The damping term is assumed to be 
non-linear and localized  to an open subset of the domain. We handle the disturbed case as an extension of  
\cite{Martinez2000}, where stability results are given with a damping term active on the full domain with no disturbances considered. We provide input-to-state types of results.
\end{abstract}
\section{Introduction}\label{intro}
Consider the damped wave equation with localized damping, with Dirichlet boundary conditions given by
\begin{equation}\label{probper}
\mathbf{(P_{dis})}\quad
\left\lbrace
\begin{array}{cccc} 
u_{tt} - \Delta u 
=-a(x)g(u_t+d)-e, &\text{in }   \mathbb{R_+} \times \Omega, \\
u=0,    &\text{on }   \mathbb{R_+} \times \partial \Omega ,\\
u(0,.)= u^0\ ,\ u_t(0,.)=u^1,
\end{array} 
\right.
\end{equation}
\ \\
where $\Omega$ is a $C^2$ bounded domain of $\mathbb{R}^2$, $d$ and $e$ stand for a damping disturbance and a globally distributed disturbance for the wave dynamics respectively. The term 
$-a(x)g(u_t+d)$ stands for the (perturbed) damping term, where $g:\mathbb{R} \longrightarrow \mathbb{R}$ is a $C^1$ non-decreasing function verifying $\xi g(\xi)>0$ for $\xi\neq 0$ while $a : \overline{\Omega} \rightarrow \mathbb{R}$ is a continuous non negative function which is bounded below by a positive constant $a_0$ on some non-empty open subset $\omega$ of 
$\Omega$. Here, $\omega$ is the region of the domain where the damping term is active, more precisely, the region where the localization function $a$ is bounded from below by $a_0$.
As for the initial condition $(u^0,u^1)$, it belongs to the standard Hilbert space 
$H^1_0(\Omega)\times L^2(\Omega)$. 

In this paper, our aim is to obtain input-to-state (ISS) type of results for 
$\mathbf{(P_{dis})}$, i.e., estimates of the norm of the state $u$  which, at once, show that trajectories tend to zero in the absence of disturbances and remain bounded by a function of the norms of the disturbances otherwise. 
\newgeometry{left=2cm,right=2cm,bottom=2cm,top=2cm}
One can refer to \cite{Sontag2008} for a thorough review of ISS results and techniques for finite dimension systems and to the recent survey 
\cite{Miron-Prieur2019} for infinite dimensional dynamical systems. 
In the case of the undisturbed dynamics, i.e., \eqref{probper} with 
$(d,e)\equiv (0,0)$, there is a vast literature regarding the stability of the 
corresponding system with respect to the origin, which is the unique 
equilibrium state of the problem. This in turn amounts to have appropriate 
assumptions on $a$ and $g$, cf. \cite{Alabau} for 
extensive references. We will however point out the main ones that we need in 
order to provide the context of our work.
To do so, we start by defining the energy of the system by
\begin{equation}\label{eq:ener}
E(t)=\frac{1}{2} \int_\Omega \left(|\nabla u|^2(t,\cdot) + u_t^2(t,\cdot) \right) \, dx ,
\end{equation}  
which defines a natural norm on the space $H^1_0(\Omega)\times L^2(\Omega)$.
Strong stabilization has been established in the early works
 \cite{Dafermos} and \cite{Haraux}, i.e., it is proved with an argument based on the Lasalle invariance principle that 
$
\lim_{t\rightarrow +\infty}E(t) = 0
$
for every initial condition in 
$H^1_0(\Omega)\times L^2(\Omega)$. However, no decay 
rate of convergence for $E$ is established since it requires in particular extra 
assumptions on $g$ and $\omega$. 

As a first working hypothesis, we will assume that 
$g'(0)>0$, classifying the present work in those that aimed at establishing results of 
exponential convergence for strong solutions.  We refer to \cite{Alabau} for the line of work where 
$g$ is assumed to be super-linear in a neighborhood of the origin (typically of 
polynomial type). Note that, in most of these works (except for the linear case) the rate of exponential decay of $E$ depends on the 
initial conditions.
That latter fact in turn relies on growth conditions of $g$ at infinity. Regarding the 
assumptions on $\omega$, they have been first put forward in the 
pioneering work \cite{Zuazua1990} on semi-linear wave equations and its 
extension in \cite{Liu1997}, where the multiplier geometric conditions 
(MGC) 
have been characterized for $\omega$ in order to achieve exponential stability. 
For linear equations, the sharpest geometrical results are obtained by 
microlocal techniques using the method of geometrical optics, cf \cite{Bardos1992} and \cite{Burq1997}. 

In this paper, our objective is to obtain results for non-linear damping terms  
and one should think of the nonlinearity $g$ not only as a mean to provide more general asymptotic behavior at infinity than a linear one but also as modeling an uncertainty of the shape of the damping term. Dealing with nonlinearities justifies  why microlocal techniques are not suited here and we will be using the multiplier method as presented e.g. in \cite{Komornik1994}.
Many results have been established in the case where $g'(0)=0$, for instance, decay rates for the energy are provided in \cite{Martinez1999} in the
localized case but the non-linearity is to have a linear growth for large
values of its arguments. Note that the
estimates as presented in \cite{Martinez1999} are not optimal in general, as for
instance in the case of a power-like growth. For general optimal energy
decay estimates and for general abstract PDEs, we refer the reader to
\cite{Alabau2005} for a general formula for explicit energy decay estimates and to
\cite{Alabau2010} for an equivalent simplified energy decay estimate with
optimality results in the finite dimensional case. However, when it comes to working under the hypothesis $g'(0)>0$,  few general results are available. One can find a rather complete presentation of the available results in 
\cite{Martinez2000}.
In particular, the proof of exponential stability along strong solutions has only been given for general nonlinearities $g$, in dimension two and in the special case of a non-localized damping with no disturbances requiring only one multiplier coupled with a judicious use of Gagliardo-Nirenberg's inequality. Our results generalize this finding in the absence of disturbances (even though it has been mentioned in \cite{Martinez2000} with no proof 
that this is the case). It has also to be noted that similar results are provided in \cite{Martinez1999} in the localized case but the nonlinearity is lower bounded by a linear function for large values of its arguments. That simplifies considerably some computations. Recall also that the purpose of \cite{Martinez1999} is instead to address issues when $g'(0)=0$ and to obtain accurate decay rates for $E$. 

Hence a possible interest of the present paper is the fact that it handles nonlinearities $g$ so that $g(v)/v$ tends to zero as $\vert v\vert$ tends to infinity with a linear behavior in a neighborhood of the origin.

As for ISS purposes, this paper can be seen as an extension to the infinite dimensional context of \cite{LiuC1995} where the nonlinearity is of the saturation type. Moreover, the present work extends to the dimension two the works \cite{Marx2017} and \cite{Marx2019}, where this type of issues have been addressed by building appropriate Lyapunov functions and by providing results in dimension one. Here, we are not able to construct Lyapunov functions and we rely instead on energy estimates based on the multiplier method, showing how these estimates change when adding the two disturbances $d$ and $e$. To develop that strategy, we must impose additional assumptions on $g'$, still handling saturation functions. As a final remark, we must recall that \cite{Martinez2000} contains other stability results in two directions. On one hand, $g'$ can simply admit a (possibly) negative lower bound and  on the other hand, the space dimension $N$ can be larger than $2$, at the price of more restrictive assumptions on $g$, in particular, by assuming quasi-linear lower bounds for its asymptotic behavior at infinity. One can readily extend the results of the present paper in both directions by eventually 
adding growth conditions on $g$. 

\section{Statement of the problem and main result}\label{statm}

In this section, we provide assumptions on the data needed to precisely define \eqref{probper}. We henceforth refer to  \eqref{probper} as the disturbed problem
$\mathbf{(P_{dis})}$. Next, we state and comment the main results of this work and discuss possible extensions.  
\\
\\
Throughout the paper, the domain $\Omega$ is a bounded open subset of $\mathbb{R}^2$ of class  $C^2$,
the assumptions on $g$ are the following.
\\ \\
$\mathbf{(H_1)}$: 
The function $g : \mathbb{R} \longrightarrow \mathbb{R}$
 is a $C^1$ non-decreasing function such that 
\begin{equation}\label{hypg}
 g(0)=0,\ \ g'(0)>0 , \ g(x)x> 0 \ \hbox{ for }x\neq 0, \ 
\end{equation}
\begin{align}\label{hypgq}
 \exists \ C>0, \ \exists \ 1<q<5,\  \forall \ |x|\geq 1,  \ |g(x)| \leq C|x|^q,
\end{align}
\begin{align}\label{hypg'}
 \exists\ C>0 ,\ \exists \ 0<m<4 ,\ \forall |x|>1,  \ \ |g'(x)| \leq C|x|^m .
\end{align}
$\mathbf{(H_2)}$: The localization function $ a : \overline{\Omega} \rightarrow \mathbb{R}\ $ is a continuous function such that 
\begin{equation} \label{a0}
a\geq 0\ \ on\ \ \Omega\ \ \text{ and }\ \exists \ a_0 >0, \ a \geq a_0\ \ on\ \ \omega.
\end{equation}
In order to prove the stability of solutions, we impose a multiplier geometrical condition (MGC) on $\omega$. It is given by the following hypothesis.
\\ \\
$\mathbf{(H_3)}$: There exists an observation point $x_0 \in \mathbb{R}^2$ for which $\omega$ contains the intersection of $\Omega$ with an $\epsilon$-neighborhood of  
\begin{equation}\label{gammadef}
\Gamma(x_0)=\{x \in \partial \Omega,\  (x-x_0).\nu (x) \geq 0\},
\end{equation}
where $\nu$ is the unit outward normal vector for $\partial \Omega$ and an $\epsilon$-neighborhood of $\Gamma(x_0)$ is defined by 
\begin{align}\label{geosect}
\mathcal{N}_\epsilon(\Gamma(x_0))= \lbrace x \in \mathbb{R}^2 \ :\ \hbox{dist}(x,\Gamma(x_0)) \leq \epsilon \rbrace.
\end{align}
Regarding the disturbances $d$ and $e$, we make the following assumptions.
\\ \\
$\mathbf{(H_4)}$: the disturbance function $d: \mathbb{R_+} \times \Omega\longrightarrow \mathbb{R}$ belongs to 
$ L^1(\mathbb{R_+},L^2(\Omega))$ and satisfies the following:	
\begin{equation}\label{lip1}
d(t,\cdot) \in  H^1_0(\Omega)\cap L^{2q}(\Omega),\ \forall t \in \mathbb{R_+},\ \ \ \   t \mapsto \int_0^t \Delta d(s,\cdot)\, ds - {d}_t(t,\cdot) \in Lip\left(\mathbb{R}_+ , H^1_0(\Omega) \right),
\end{equation}
where $Lip$ denotes the space of Lipschitz continuous functions. We also
impose that the following quantities
\begin{align}\label{bornc1}
&C_1(d) = \int_0^\infty \int_\Omega (|d|^2+|d|^{2q}) \, dx\, dt  ,\ \  \ \  
C_2(d)=	\int_0^\infty 
\int_{\Omega} |d|^m\, (d_t)^2 \, dx\, dt,\ \  \ \ 
 \notag
\\ 
&C_3(d)=\int_0^\infty  \int_{\Omega} (d_t)^2 \, dx\,dt,\ \  \ \ 
C_4(d)=\int_0^\infty \left( \int_{\Omega}  |{d}_t|^{2\left( \frac{p}{p-1} \right)} \, dx \right)^{\left( \frac{p-1}{p} \right)} dt ,
\end{align}
are all finite, where $p$ is a fixed real number so that, if\   $ 0<m\leq 2$,\  then\  $p>\frac{2}{m}$
and if $2<m<4$, then  $p\in (1,\frac{m}{m-2})$.
\begin{Remark}
	The fact that $d$ belongs to
	$ L^1(\mathbb{R_+},L^2(\Omega))$ means that the following quantity is finite
	\begin{align}
	&C_5(d) = \int_0^\infty  || d||_{L^2(\Omega)}\, dt,
	\end{align}
	which implies that the following quantity is also finite
	\begin{align}
	 C_6(d) = \int_0^\infty  \int_\Omega |d| \, dx\, dt.
	\end{align}
\end{Remark}
$\mathbf{(H_5)}$: The disturbance function $e:  \mathbb{R_+} \times \Omega\longrightarrow \mathbb{R}$ belongs to $W^{1,1}(\mathbb{R_+},L^2(\Omega))$ and satisfies the following	
\begin{equation}\label{lip2}
e \ \in \   Lip\left(\mathbb{R}_+ , H^1_0(\Omega)\right), \ \ \ e(0,.)\in L^2(\Omega), \ \ \ \bar{C}_1(e) = \int_0^\infty  \int_{\Omega} e^2 \, dx\, dt < \infty.
\end{equation} 
\begin{Remark}
	The fact that $e$ belongs to
	$ W^{1,1}(\mathbb{R_+},L^2(\Omega))$ means that the following quantities are finite
	\begin{align}
	&\bar{C}_2(e) = \int_0^\infty  || e(t,\cdot)||_{L^2(\Omega)}\, dt, \ \ \ \ \bar{C}_3(e) = \int_0^\infty  || e_t(t,\cdot)||_{L^2(\Omega)}\, dt.
	\end{align}
\end{Remark}

\begin{Remark}
In the rest of the paper, we will use various symbols $C$, $C_u$ and  
$C_{d,e}$ which are constants independent of the time $t$. 
However, it is important to 
stress that these symbols have specific dependence on other parameters of the 
problem. More precisely, the symbol $C$ will be used to denote positive constants 
independent of initial conditions and disturbances, i.e., only depending on the
domains $\Omega,\omega$ and the functions $a$ and $g$. The symbol $C_u$ 
denotes a generic $\mathcal{K}$-function of the norms of the initial 
condition $(u_0,u_1)$ and  similarly the symbol $C_{d,e}$
denotes a generic $\mathcal{K}$-function of the several quantities $C_i(d)$ and $\bar{C}_i(e)$. Here $\mathcal{K}$ denotes the set of continuous increasing functions $ \gamma :\ 
\mathbb{R}_{+}\to \mathbb{R}_{+}$ with $ \gamma (0)=0$, cf. \cite{Miron-Prieur2019}.
\\
Moreover, in the course of intermediate computations, we will try to keep all the 
previous constants as explicit as possible in terms of the norms of 
the initial condition and the $C_i(d)$ and $\bar{C}_i(e)$ in order to keep track of the nature of 
generic constants. We will use the latter generic mainly in the statements of the results.
\end{Remark}
Before we state the main results, we define the notion of a strong solution of $\mathbf{(P_{dis})}$. To do so, we start by giving an equivalent form of $\mathbf{(P_{dis})}$ : 
	\\
	\\
Define for every $(t,x)\in \mathbb{R_+}\times\Omega$,  $\bar{d}(t,x)=\int_0^t d(s,x) 
ds$. We translate $u$ in $\mathbf{(P_{dis})}$ as $v=u+\bar{d}$, it is immediate to see 
that $\mathbf{(P_{dis})}$ is equivalent to
the following problem:
\begin{equation}\label{probpert}
\left\lbrace
\begin{array}{cccc} 
v_{tt} - \Delta v + a(x)g(v_t)=\tilde{e},
&\text{in }   \mathbb{R_+} \times \Omega, \\
v=0,    &\text{on }  \mathbb{R_+} \times \partial \Omega,\\
v(0,.)= v^0,\ \  v_t(0,.)=\bar{u}^1,
\end{array} 
\right.
\end{equation}
where $\tilde{e}= d_t -\Delta \bar{d}-e$, $v^0=u^0$ and $v^1=u^1 +d(0,.)$.
\\
Define the unbounded operator
\begin{align}\label{Adef}
A : \ H=H^1_0(\Omega)\times L^2(\Omega) &\longrightarrow H , \notag\\
(x_1,x_2) &\longmapsto  (x_2,-\Delta x_1 +ag(x_2)), 
\end{align}
with domain
\begin{align*}
D(A) = \left(H^2(\Omega)\cap H^1_0(\Omega)\right) \times H^1_0(\Omega).
\end{align*}
For $t\geq 0$, set 
$$
U(t)=\begin{pmatrix}u(t,\cdot) \\\ u_t(t,\cdot)\end{pmatrix},\quad
V(t)=\begin{pmatrix}v(t,\cdot) \\\ v_t(t,\cdot)\end{pmatrix},\quad
D(t)=\begin{pmatrix}\bar{d}(t,\cdot) \\\ d(t,\cdot)\end{pmatrix},\quad
G(t) = \begin{pmatrix} 0 \\ \tilde{e} \end{pmatrix}.
$$
Notice that $G \in Lip(\mathbb{R}_+, L^2(\Omega)\times H^1_0(\Omega))$.
Then 
Problem \eqref{probpert} can be written as  
\begin{equation} \label{veq}
V_t(t) = AV(t) + G(t), 
\ \ V(0)=V_0\ = \begin{pmatrix} v^0 \\ v^1 \end{pmatrix}.
\end{equation}
A strong solution of \eqref{veq} in the sens of \cite{Brezismmo} is a function $V\in 
C(\mathbb{R_+},H)$, absolutely continuous in every compact of $\mathbb{R_+}$, 
satisfying $V(t) \in D(A), \forall t\in \mathbb{R_+}$ and satisfying \eqref{veq} almost 
everywhere in $\mathbb{R_+}$. On the other hand, the hypotheses satisfied by $d$ 
imply that $D(t) \in D(A)$ for every $t\in\mathbb{R}_+$. Since $U=V-D$, we can now 
give the following definition for a strong solution of $\mathbf{(P_{dis})}$.

\begin{Definition} (Strong solution of $\mathbf{(P_{dis})}$.)\\
	A strong solution $u$ of  $\mathbf{(P_{dis})}$ is a function $u\in C^1(\mathbb{R}_+, L^2(\Omega)) \cap C(\mathbb{R}_+,H^1_0(\Omega)) $ such that
	$t \mapsto u_t(t,\cdot)$ is absolutely continous in every compact of $\mathbb{R_+}$.  
For all $t\in \mathbb{R_+},$ $(u(t,\cdot),u_t(t,\cdot)) \in D(A)$ and  $u(t,\cdot)$ satisfies $\mathbf{(P_{dis})}$ for almost all $t\in \mathbb{R_+}$.
\end{Definition}
\ \\
We gather our findings in the following theorem regarding the disturbed system $\mathbf{(P_{dis})}$.



\begin{Theorem} \label{thrmdis}
Suppose that Hypotheses $\mathbf{(H_1)}$ to  
$\mathbf{(H_5)}$ are satisfied. Then, given 
$(u^0,u^1) \in (H^2(\Omega)\cap H^1_0(\Omega)) \times H^1_0(\Omega) $, Problem $\mathbf{(P_{dis})}$ has a unique strong solution $u$.
Furthermore, the following energy estimate holds:
{\begin{align} \label{gt1}
E(t) \leq  (C+C_u)E(0) e^{-\frac{t-1}{C_u+C}}  + C_{d,e}(C_u+1),
\end{align}}
where the positive constant $C_u$ depends only on the initial conditions  
and the positive constant $C_{d,e}$ depends only on the disturbances $d$ and $e$.
\end{Theorem}
\begin{Remark}\textbf{(Comments and extensions)}
\begin{itemize}
\item Theorem~\ref{thrmdis} holds true if  the Lipschitz assumptions in \eqref{lip1} and \eqref{lip2} are replaced by bounded variation ones.
\item  In the case where the disturbances are both zero ($d\equiv 0$ and $e\equiv0$), 
Theorem~\ref{thrmdis} holds without the hypothesis on $g'$ given by \eqref{hypg'} (i.e. 
no restriction on $q$ in \eqref{hypgq}) and the hypothesis given by \eqref{hypgq} can be then weakened to the following one
\begin{align*}
\exists\ C>0, \ \exists\  q>1,\ \forall \ |x|\geq 1, \ |g(x)|\leq C|x|^q.
\end{align*}
It is clear that if $g$ satisfies the last part of the condition above for $0\leq q\leq 1$, it would still satisfy it for any $q>1$.
\item The geometrical condition MGC imposed in $\mathbf{(H_3)}$ can be readily reduced to the weaker and more general  MGC introduced in \cite{Liu1997} and called piecewise MGC in \cite{Alabau}.

\end{itemize}
\end{Remark}
\begin{Remark}
Note that \eqref{gt1} is an ISS-type estimate but it fails to be a strict one (let say 
in the sense of Definition $1.6$ in \cite{Miron-Prieur2019}) for two facts. First of all, the estimated 
quantity $E$ is the norm of a trajectory in the space $H^1_0(\Omega)\times 
L^2(\Omega)$ while the constant $C_u$ depends on the initial condition by
its norm in the smaller space $(H^2(\Omega)\cap H^1_0(\Omega)) \times 
H^1_0(\Omega)$. This difference seems unavoidable since in the undisturbed case 
exponential decay can be proved only for strong solutions as soon as the nonlinearity $g$ is not assumed to be bounded below at infinity by a linear function. As a matter of fact, it would be interesting to prove that strong stability is the best convergence result one could get for weak solutions, let say with damping functions $g$ of saturation type functions and in dimension at least two.
\\
The second difference lies in the second term in \eqref{gt1}, 
namely it is not just a $\mathcal{K}$-function of the norms of the disturbances. We 
can get such a result if we have an extra assumption on $g$, typically $g$ of growth 
at most linear at infinity (i.e., $q=1$) with bounded derivative (i.e., $m=0$). In 
particular, this covers the case of regular saturation functions (increasing bounded 
functions $g$ with bounded derivatives).
\end{Remark}
We give now the proof of the well-posedness part of Theorem \eqref{thrmdis}.
\\ \\
\textbf{\emph{Proof of the well-posedness:}} The argument is standard since $-A$, 
where $A$ is defined in \eqref{Adef}, is a maximal monotone operator on $ 
H^1_0(\Omega)\times L^2(\Omega)$ (cf. for instance \cite{harauxevo} for a proof).
We can apply Theorem~3.4 combined with Propositions~3.2 and Propositions~3.3 in 
\cite{Brezismmo} to \eqref{veq}, which immediately proves the results of the well-
posedness part.
\begin{flushright}
	\begin{small}
$\blacksquare$
\end{small}
\end{flushright}
\begin{Remark}
	In \cite{Martinez2000}, the domain of the operator has been chosen as  
	\begin{align*}
 \lbrace (u,v) \in H^1_0(\Omega) \times H^1_0(\Omega) \ : \ -\Delta u +g(v) \in L^2(\Omega) \rbrace.
	\end{align*}
	However, in dimension two, taking the domain of $A$ in the case where $d=e=0$
	as $Z=\lbrace (u,v) \in H^1_0(\Omega) \times H^1_0(\Omega):-\Delta u +a(x)g(v) \in L^2(\Omega) \rbrace$ or as $(H^2(\Omega) \cap H^1_0(\Omega)) \times H^1_0(\Omega) $ is equivalent.	Indeed, using the hypothesis given by \eqref{hypgq}, we have that $|g(v)| \leq C|v|^q $ for $|v|<1$, which means when combining it with the fact that $g(0)=0$ that $|g(v)| \leq C|v|^q + C|v| $ for all $v$. From Gagliardo-Nirenberg theorem (see in Appendix) we have for $v\in H^1_0(\Omega)$ that
	\begin{align*}
	\Vert  v\Vert _{L^{2q}(\Omega)}^{2q} \leq C \Vert  v \Vert _{H^1_0(\Omega)}^{2q-2} \Vert  v \Vert _{L^{2}(\Omega)}^2,
	\end{align*}
	which means that
	\begin{align*}
	{	\Vert  g(v) \Vert _{L^2(\Omega)}^2} &= \int_\Omega |g(v)|^2 dx \leq C\int_\Omega \left(|v|^q + |v|\right)^{2} dx \leq  C\Vert  v \Vert _{L^2q(\Omega)}^{2q}+C\Vert  v \Vert _{L^2(\Omega)}^2 
	\\
	&\leq  \Vert  v \Vert _{H^1_0(\Omega)}^{2q-2} \Vert  v \Vert _{L^{2}(\Omega)}^2+C\Vert  v \Vert _{L^2(\Omega)}^2 <+\infty\ \ (\text{since}\ v\in H^1_0(\Omega) ),
	\end{align*}
	i.e., 
$g(v) \in L^2(\Omega)$. Then, by using Lemma 
	\ref{lemma2per} (with $(d,e)\equiv (0,0)$), we have that $-\Delta u + ag(v) \in L^2(\Omega)$, which 
	means that $\Delta u \in L^2(\Omega)$. On the other hand, $\Vert  \Delta u 
	\Vert _{L^2(\Omega)}$ is an equivalent norm to the norm of $H^2(\Omega)\cap H^1_0(\Omega)$ and $\Omega$ is of class $C^2$ (the proof is a 
	direct result of Theorem 4 of Section 6.3 in \cite{Evans}).
	We can finally conclude that $Z$ is nothing else but $(H^2(\Omega) \cap H^1_0(\Omega)) \times H^1_0(\Omega) $. 
\end{Remark}
\section{Proof of the energy estimate \eqref{gt1} }\label{proofs} \label{mainresult}

 To prove the energy estimate given by \eqref{gt1}, we are going to use the multiplier method combined with a Gronwall lemma and other technical lemmas given in this section. We will be referring to \cite{Martinez1999} and \cite{Martinez2000} in several computations since our problem is a generalization of their strategy to the case where the disturbances $(d,e)$ are present. 
\\ \\
We start with the following lemma stating that the energy $E$ is bounded along trajectories of $\mathbf{(P_{dis})}$.
\begin{Lemma} \label{lemma1per} 
Under the hypotheses of Theorem \eqref{thrmdis}, the energy of a strong solution of Problem $\mathbf{(P_{dis})}$, satisfies 
\begin{equation}\label{nonincper}
E'(t)= - \int_\Omega a u_tg(u_t+d)\, dx - \int_\Omega u_t e \, dx,\quad \forall t\geq 0.
\end{equation}
Furthermore, there exist positive constants $C$ and $C_{d,e}$ such that
\begin{equation}\label{boundE'}
E(T) \leq CE(S)+C_{d,e},\quad \forall\  0\leq S\leq T.
\end{equation}
\end{Lemma}   
\textbf{\emph{Proof of Lemma~\ref{lemma1per}: }}
Equation \eqref{nonincper} follows after multiplying the first equation of \eqref{probper} by $u_t$ and performing standard computations. Notice that we do not have the dissipation of $E$ since the sign of $E'$ is not necessarily constant.
To achieve \eqref{boundE'}, we first write
\begin{align}\label{bnd0}
- \int_\Omega a u_tg(u_t+d)\, dx =- \int_{|u_t|\leq |d|} a u_tg(u_t+d) \, dx 
-  \int_{|u_t|> |d|} a u_tg(u_t+d)\, dx.
\end{align}
On one hand, from \eqref{hypg} and the fact that $(u_t+d)$ and $u_t$ have the same sign if 
$|u_t|> |d|$, we deduce that
\begin{align}\label{bnd1}
-\int_{|u_t|> |d|} au_tg(u_t+d) \, dx \leq 0.
\end{align}
On the other hand, since $g$ is non-decreasing, has linear growth in a neighborhood of zero by \eqref{hypg}, and satisfies \eqref{hypgq}, it follows that
\begin{align}\label{bnd2}
 -\int_{|u_t|\leq |d|} a u_tg(u_t+d) \, dx &\leq C \int_{|u_t|\leq |d|} |d||g(|2d|)| \, dx \leq C \int_\Omega |d||g(|2d|)| \, dx \notag
 \\
& \leq C \int_{|d|<1}  |d||g(2d)| dx +C\int_{|d|\geq 1} |d||g(2d)| dx 
\notag 
\\
&\leq C \int_{|d|<1}  |d|^2 dx + C\int_{|d|\geq 1} |d|^{q+1} dx  
\notag\\
&\leq C  \int_\Omega (\vert d\vert^2+|d|^{2q}) \, dx.
\end{align}
Combining \eqref{bnd0}, \eqref{bnd1}, \eqref{bnd2} and \eqref{nonincper}, we obtain that
\begin{align}\label{bnd0'}
E' \leq C \int_\Omega (\vert d\vert^2+|d|^{2q}) \, dx -  \int_{\Omega} u_t e \, dx\, dt.
\end{align}
Using Cauchy-Schwarz inequality,
\begin{align*}
E'&\leq C\int_\Omega (|d|^2+|d|^{2q}) \, dx  + \left( \int_{\Omega} |e |^2  \, dx\, dt \right)^\frac{1}{2} \left( \int_{\Omega} |u_t |^2 \, dx\right)^\frac{1}{2}\\
&\leq C\int_\Omega (|d|^2+|d|^{2q}) \, dx  +C \Vert  e \Vert _{L^2(\Omega)} \sqrt{E},
\end{align*}
then integrating between two arbitrary non negative times $S\leq T$, we get
\begin{small}
	\begin{align*}\label{bnd0''}
 E(T) \leq E(S) +CC_1(d) + C\int_S^T \Vert  e \Vert _{L^2(\Omega)} \sqrt{E} dt,
\end{align*}
\end{small}
which allows us to apply Theorem \ref{gronw2lem}
and conclude that
\begin{align*}
E(T) \leq  CE(S)+CC_1(d)+C\bar{C}_2(e)^2 = CE(S)+C_{d,e}.
\end{align*}
Hence, the proof of Lemma \ref{lemma1per} is completed.
\begin{flushright}
\small{$\blacksquare$}
\end{flushright}
\begin{Remark}
	In the absence of disturbances, in other words when $d=e=0$ we have that: 
	\begin{equation}
	E'(t)= - \int_\Omega a u_tg(u_t)\, dx,\quad \forall t\geq 0,
	\end{equation}
	and thus the energy $E$ is non increasing by using \eqref{hypg}. That latter fact simplifies the proof of exponential decrease in this case.
\end{Remark}
We provide now an extension of Lemma 2 in \cite{Martinez2000} to the context of $\mathbf{(P_{dis})}$.
\begin{Lemma} \label{lemma2per}
Under the hypotheses of Theorem \ref{thrmdis}, for every solution of Problem 
$\mathbf{(P_{dis})}$ with initial conditions $(u^0,u^1)\in (H^2(\Omega)\cap 
H^1_0(\Omega))\times H^1_0(\Omega)$, there exist explicit positive constants
$C_u$ and $C_{d,e}$ such that 
\begin{equation} \label{Clemmaper}
\forall t \geq 0,\ \ \ \Vert  - \Delta u(t,\cdot)+a(\cdot)g(u_t(t,\cdot) + d(t,\cdot)) + e(t,\cdot)\Vert  ^2_{L^2(\Omega)} + \Vert  u_t(t,\cdot) \Vert  ^2_{H^1_0(\Omega)} \leq   C_u+C_{d,e}.
\end{equation}
\end{Lemma}
\textbf{\emph{Proof of Lemma~\ref{lemma2per}:}} We set  $w:=u_t$, where $u$ is the strong solution of $\mathbf{(P_{dis})}$. We know that $w(t) \in H^1_0(\Omega)$ for every $t\geq 0$. Moreover, it is standard to show that  
$w(t)$ satisfies in the distributional sense the following problem:
\begin{equation}\label{wdiffper}
\left\lbrace
\begin{array}{cccc} 
{w}_{tt}-\Delta w + ag'(w+d)(w_t +{d}_t)+{e}_t=0, &\text{in } \Omega \times \mathbb{R_+}, \\
w=0,    &\text{on } \partial \Omega \times \mathbb{R}_+ ,\\
w(0)=u^1 , \ \ w_t(0)=\Delta u^0 -g(u^1+d(0))-e(0).
\end{array} 
\right.
\end{equation}
Set $E_w(t)$ to be the energy of $w$ for all $t\geq 0$. It is given by 
\begin{equation*}
E_w(t)=\frac{1}{2} \int_{\Omega} (w_t^2(t,x) +|\nabla w(t,x) |^2) \, dx.
\end{equation*}
Using $w_t$ as a test function in \eqref{wdiffper}, then performing standard computations,
we derive
\begin{align}\label{eq:perEw1}
E_w(t)-E_w(0)=&-\int_0^t \int_{\Omega} (a g'(w+d)({d}_t+w_t)w_t+ {e}_tw_t) \, dxd\tau.
\end{align}
Let ${I}:=\int_0^t \int_{\Omega} a(.) g'(w+d)({d}_t+w_t)w_t\, dxd\tau$. We split the domain $\Omega$ in $I$ according to whether $|{d}_t|\leq |w_t|$ or not. Clearly the part corresponding to 
$|{d}_t|\leq |w_t|$ is non negative since $g'\geq 0$, $a\geq 0$ and $({d}_t+w_t)$ and $w_t$ have the same sign. 
From \eqref{hypg'}, one has the immediate estimate
$$
g'(a+b)\leq C(1+\vert a+b\vert^m)\leq C(1+\vert a\vert^m+\vert b\vert^m),
\quad \forall a,b\in \mathbb{R}.
$$
Using the above, we can rewrite \eqref{eq:perEw1}
as
\begin{align}
E_w(t)-E_w(0)&\leq \int_0^t \int_{|{d}_t|>|w_t|} a g'(w+d)({d}_t+w_t)w_t\, dxd\tau +\int_0^t \int_{\Omega} |{e}_t||w_t| \, dx d\tau  \notag \\
&\leq C 
\int_0^t \int_{\Omega} g'(w+d){d}_t^2 \, dxd\tau +C\int_0^t ||{e}_t||_{L^2(\Omega)} \sqrt{E_w} d\tau
\notag \\
&\leq C \int_0^t 
\int_{\Omega} (1+|w|^m+|d|^m){d}_t^2 \, dxd\tau +C\int_0^t \Vert{e}_t\Vert_{L^2(\Omega)} \sqrt{E_w} d\tau.\label{ww0}
\end{align}
Using Hölder's inequality,  
\begin{align}\label{ww1}
\int_0^t \int_{\Omega} |w|^m {d}_t^2 \, dxd\tau \leq \int_0^t \left( \int_{\Omega} |w|^{pm} \, dx \right)^\frac{1}{p}\left( \int_{\Omega} |{d}_t|^{2p'} \, dx \right)^\frac{1}{p'} d\tau,
\end{align}
with $p$ defined in \eqref{bornc1} and $p'> 1$ is its conjugate exponent given 
by $\frac{1}{p}+\frac{1}{p'}=1$. Thanks to the assumptions on $p$, one can use Gagliardo-Nirenberg's inequality for $w$ to get
\begin{align} \label{gagniw}
\left(\int_{\Omega} |w(t,x)|^{pm} \, dx \right)^\frac{1}{p} \leq C E_w(t)^{\frac{m \theta}{2}} E(t)^{\frac{(1-\theta)m}{2}},\quad t\geq 0, 
\end{align}
where $\theta =1- \frac{2}{mp}$.
Combining \eqref{gagniw}, \eqref{ww1} and \eqref{ww0}, it follows that
\begin{align}\label{ww3}
&E_w(t)-E_w(0) \leq  C \int_0^t  E_w^{\frac{m \theta}{2}} E^{\frac{(1-\theta)m}{2}} \int_{\Omega} \left( |{d}_t|^{2p'} \, dx \right)^\frac{1}{p'} d\tau\notag\\
& \ \ \ \ + \int_0^t 
\int_{\Omega} (1+|d|^m) {d}_t^2 \, dx d\tau 
+C\int_0^t ||{e}_t||_{L^2(\Omega)} \sqrt{E_w} d\tau.
\end{align}
Note that $\frac{m \theta}{2} <1$. Setting $h_1(t)=\int_{\Omega} \left( |{d}_t|^{2p'} \, dx \right)^\frac{1}{p'} $, $h_2(t) = ||{e}_t||_{L^2(\Omega)}  $ and using \eqref{boundE'},
\eqref{ww3} becomes 
\begin{align}\label{ww3'}
 E_w(t)  \leq E_w(0) +C_2(d)+C_3(d) + \left(C_u+ C_{d,e}\right) & \int_0^t E_w^{\frac{m \theta}{2}} h_1(s) ds  +C\int_0^t h_2(s) \sqrt{E_w} ds.
\end{align}
We know that 
\begin{align}
\int_0^\infty h_1(t)\, dt = C_4(d) <\infty, \ \ \ \int_0^\infty h_2(t)\, dt = \bar{C}_3(e)<\infty.
\end{align}
We can now apply Theorem \ref{gronw2lem} on \eqref{ww3'} with 
\begin{align*}
S=0,\  T=t, \ \alpha_1 = \frac{m\theta}{2},\ & \alpha_1 = \frac{1}{2},\
F(\cdot)=E_w(\cdot),\ C_3=C_2(d)+C_3(d),\ C_1 = C_u+ C_{d,e}, \ C_2=C.
\end{align*}
We obtain the following bound for $E_w(\cdot)$:
\begin{align}\label{ew}
E_w(t) \leq \max\Big(2(E_w(0)+C_2(d)+C_3(d)),(2\tilde{C})^{\frac{1}{1-\alpha}}\Big),
\end{align}
 where $\tilde{C}:=C_1\Vert h_1\Vert_1+C_2\Vert h_2\Vert_1$ and $\alpha:=\max(\alpha_1,\alpha_2)$ if $2\tilde{C}\geq 1$ or 
$\alpha:=\min(\alpha_1,\alpha_2)$ if $2\tilde{C}< 1$. \\
\\
It is clear that $\tilde{C}=(C_u+C_{d,e})C_4(d)+C\bar{C}_3(e)\leq C_u+C_{d,e}$.
One then rewrites \eqref{ew} as 
\begin{align}\label{eww}
E_w(t) \leq 2(E_w(0)+C_2(d)+C_3(d))+(C_u+C_{d,e})^{\frac{1}{1-\alpha}}.
\end{align}
Note that for $t \geq 0$ one obviously has that
\begin{align*}
	E_w(t)&=  \frac{1}{2} \int_{\Omega} (w_t^2(t,x) +|\nabla w(t,x) |^2) \, dx
	\\
	&=  \frac{1}{2} \left( || u_{tt}(t,\cdot)||^2_{L^2(\Omega)}+\Vert  u_t(t,\cdot) \Vert  ^2_{H^1_0(\Omega)}\right)
	\\
	&=\Vert  - \Delta u(t,\cdot)+a(\cdot)g(u_t(t,\cdot) + d(t,\cdot)) + e(t,\cdot)\Vert  ^2_{L^2(\Omega)} + \Vert  u_t(t,\cdot) \Vert  ^2_{H^1_0(\Omega)}.
\end{align*}
The conclusion of the lemma follows since, by taking into account \eqref{hypgq}, it is 
clear that $E_w(0)\leq C_u+C_{d,e}$.
%
%
\begin{flushright}
\small{$\blacksquare$}
\end{flushright}
We next provide the following important estimate based on Gagliardo-Nirenberg theorem:
\begin{Lemma}\label{gagnircorper}
For all $q> 2 $, a strong solution  $u$ of $\mathbf{(P_{dis})}$ satisfies
\begin{equation}\label{gagniper}
\Vert  u_t(t,\cdot) \Vert _{L^q(\Omega)}^q \leq  (C_u + C_{d,e}) E(t),\quad t\geq 0.
\end{equation}
\end{Lemma}
\textbf{\emph{Proof of Lemma~\ref{gagnircorper}: } }
We derive immediately from \eqref{Clemmaper} that
$\Vert  u_t \Vert _{H^1_0(\Omega)} \leq C_u + C_{d,e}$. Then, using Gagliardo-Nirenberg's theorem, it follows that, for every $t\geq 0$,
\begin{align}\label{twoper}
\Vert  u_t(t,\cdot)\Vert _{L^q(\Omega)}^q \leq C \Vert  u_t(t,\cdot) \Vert _{H^1_0(\Omega)}^{q-2} \Vert  u_t(t,\cdot) \Vert _{L^2(\Omega)}^2 \leq  (C_u + C_{d,e}) E(t). 
\end{align}
\begin{flushright}
\small{$\blacksquare$}
\end{flushright}
We have all the tools now to start the proof of the second part of Theorem~\ref{thrmdis}. The stability result will be achieved as a direct consequence of the following proposition:
\begin{Proposition}\label{propper}
Suppose that the hypotheses of Theorem \eqref{thrmdis} are satisfied,
then the energy $E$ of the strong solution $u$ of $\mathbf{(P_{dis})}$ with $(u^0,u^1) \in (H^2(\Omega)\cap H^1_0(\Omega)) \times H^1_0(\Omega)) $, satisfies the following estimate: 
{\begin{equation} \label{lem73per}
\int_S^T E(t)\, dt \leq (C_u+C) E(S) +(1+C_u)C_{d,e},
\end{equation}}
where the positive constant $C_u$ depends only on the initial condition, the positive constant $C_{d,e}$ depends only on the disturbances $d$ and $e$ respectively and $C$ is a positive real constant. 
\end{Proposition}

\subsubsection{Proof of Proposition~\ref{propper}}
We now embark on an argument for Proposition~\ref{propper}. It is based on the use of several multipliers that we will apply to the partial differential equation of \eqref{probper}. For that purpose, we need to define several functions associated with $\Omega$. 
\\ \\
Let $(u^0,u^1) \in (H^2(\Omega)\cap H^1_0(\Omega))\times H^1_0(\Omega) $,
$S\leq T$ two non negative times and $x_0 \in \mathbb{R}^2$ an observation point. Define $\epsilon_0$, $\epsilon_1$ and $\epsilon_2$ three positive real constants such that $\epsilon_0 < \epsilon_1 < \epsilon_2 <\epsilon$ where $\epsilon$ is the same defined in \ref{geosect}. Using $\epsilon_i$, we define $Q_i$ for $i=0,1,2$ as $
Q_i = \mathcal{N}_{\epsilon_i}[\Gamma(x_0)]$. 
\\[3pt]
Since $\overline{(\Omega \setminus Q_1)} \cap \overline{Q_0}=\emptyset $, we are allowed to define a function $\psi \in C^\infty_0(\mathbb{R}^2)$ such that 
\begin{equation*}\label{psi}
\begin{cases}  
0\leq \psi \leq 1 ,
\\ 
\psi = 1\ \hbox{on}\ \bar{\Omega}\setminus Q_1 ,
\\ 
\psi=0\ \hbox{on}\ Q_0.
\end{cases}
\end{equation*}
We also define the $C^1$ vector field $h$ on $\Omega$ by
\begin{equation}\label{hdef}
h(x):= \psi (x)(x-x_0).
\end{equation}
When the context is clear, we will omit the arguments of $h$.
\\
\\
We use the multiplier $M(u):=h\nabla u+\frac{u}2$ to deduce the following first estimate:
\begin{Lemma}\label{lemcombper}
Under the hypotheses of Proposition \ref{propper}, we have the following inequality:
\begin{align}\label{ineq2per}
\int_S^T E\ dt \leq   \underbrace{\left|\left[ \int_{\Omega} u_tM(u)\, dx \right]_S^T\right| }_{\mathbf{T_1}} +C& \underbrace{\int_S^T \int_{\Omega\cap Q_1} |\nabla u|^2  \, dx\, dt }_{\mathbf{T_2}}
+  \underbrace{ \left| \int_S^T \int_{\Omega} a g(u_t+d) M(u) \, dx\, dt \right|}_{\mathbf{T_3}} \notag
\\
&+\underbrace{ \left| \int_S^T \int_{\Omega} e M(u) \, dx\, dt\right|}_{\mathbf{T_4}} +C\underbrace{\int_S^T \int_{\omega} u_t^2\, dx\, dt}_{\mathbf{T_5}} ,
\end{align}
where $h$ is defined in \eqref{hdef} and $M(u)$ is the multiplier given by $h.\nabla u + \frac{u}{2}.
$
\end{Lemma}
\textbf{\emph{Proof of Lemma~\ref{lemcombper}. }}
The proof  is based on multiplying $\mathbf{(P_{dis})}$ by the multiplier $M(u)$ and integrating on $[S,T]\times \Omega$. Then, we follow the steps that led to the proof of equation $(3.15)$ in \cite{Martinez1999} except that we take $\sigma=0$ and $\phi(t)=t$ in the beginning and we replace $\rho(x,u_t)$ by $a(x)g(u_t+d) + e$.
\begin{flushright}
\small{$\blacksquare$}
\end{flushright}
\begin{Remark}
From now on, whenever we refer to a proof in \cite{Martinez1999}, we refer to the steps of the proof with the change of $\sigma=0$ and $\phi(t)=t$ as well as replacing $\rho(x,u_t)$ by $a(x)g(u_t+d) + e$.
\end{Remark}
The goal now is to estimate the terms ${T_1}$ to ${T_5}$. 
\begin{Lemma} \label{es1}
	Under the hypotheses of Proposition \ref{propper}, there exists a positive constant $C$ such that
	\begin{equation}\label{estim1per}
{T_1}\leq CE(S)+ C_{d,e}.
	\end{equation}
\end{Lemma}
\textbf{\emph{Proof of Lemma \ref{es1}:}}
\noindent Exactly as the proof of equation $(5.14)$ in \cite{Martinez1999} except that we use \eqref{boundE'} in the very last 
step since we do not have the non-increasing of the energy here. We obtain \eqref{estim1per}.
\begin{flushright}
	\small{$\blacksquare$}
\end{flushright}
The estimation of ${T_2}$ requires more work and it is given in the following 
lemma:
\begin{Lemma}\label{es2per}
Under the hypotheses of Proposition \ref{propper}, $T_2$ is estimated by
\begin{align}\label{ineqnablaper}
{T_2}\leq C\eta_0
\int_S^T E\ dt +\frac{C}{\eta_0}\int_S^T  \int_\omega u_t^2  \, dx\, dt +\frac{1}{\eta_0}\left(C+C_u+C_{d,e}\right) E(S) 
+ \frac{1}{\eta_0^5}\left(C_{d,e}C_u+ C_{d,e}\right),
\end{align}
where $0<\eta_0<1$ is an arbitrary real positive number to be chosen later.
\end{Lemma}
\textbf{\emph{Proof of Lemma~\ref{es2per}: }}
 The argument requires a new multiplier, namely  $\xi u$, where the  function $\xi \in C^\infty_0(\mathbb{R}^2)$ 
is defined by 
\begin{equation} \label{xidef}
\begin{cases}  
0\leq \xi \leq 1,
\\ 
\xi = 1\ \hbox{on}\ Q_1 ,
\\ 
\xi=0\ \hbox{on}\ \mathbb{R}^2\setminus Q_2.
\end{cases}
\end{equation}
Such a function $\xi$ exists since $\overline{\mathbb{R}^2\setminus 
	Q_2}\cap \overline{Q_1}=\emptyset$. Using the multiplier $\xi u$ and following the steps in the proof of Lemma~$9$ in 
\cite{Martinez1999}, yields the following identity:
\begin{small}
	\begin{align}
	\int_S^T  \int_\Omega \xi |\nabla u|^2  \, dx\, dt 
	&= \int_S^T \int_\Omega \xi |u_t|^2  \, dx\, dt+\frac{1}{2}\int_S^T \int_\Omega \Delta \xi u^2 \, dx\, dt - \left[  \int_\Omega \xi u u_t \, dx  \right]_S^T \notag
	\\
	&\ \ \ \ \ \ -\int_S^T \int_\Omega \xi u \left[a(x)g(u_t+d) + e\right]\  \, dx\, dt. \label{id3per}
	\end{align}
\end{small}
Combining the fact that $\Delta \xi$ is bounded and the definition of $\xi$, we derive from \eqref{id3per} that
\begin{align}
{T_2} \leq \int_S^T  \int_{\Omega\cap Q_2} |u_t|^2  \, dx\, dt+ 
\underbrace{\left| \left[  \int_{\Omega\cap Q_2} u u_t \, dx \right]_S^T \right|}_{{S_1}}+C\underbrace{\int_S^T \int_{\Omega\cap Q_2} u^2 \, dx\, dt}_{{S_2}}
\notag
\\
+\underbrace{\int_S^T  \int_{\Omega} |u ag(u_t+d)|  \, dx\, dt}_{{S_3}} + 
\underbrace{\int_S^T \int_\Omega |u e| \, dx\, dt}_{{S_4}}. \label{Aper}
\end{align}
First, note that the first term of \eqref{Aper} is upper bounded by 
$\int_S^T  \int_{\omega} |u_t|^2  \, dx\, dt$ since $\Omega\cap Q_2\subset\omega$. Left to estimate the other terms in the right-hand side of \eqref{Aper}. We start by treating ${S_1}$. We easily get the following estimate by using Young and Poincar\'e inequalities:
\begin{align}
\int_{\Omega\cap Q_2} |u u_t| \, dx \leq \frac{1}{2}\int_{\Omega\cap Q_2} |u|^2 \, dx + \frac{1}{2}\int_{\Omega\cap Q_2} |u_t|^2 \, dx \leq C E. \label{AAA}
\end{align}
Using \eqref{boundE'} with \eqref{AAA} we obtain the estimation of $S_1$ given by
\begin{align}\label{xiiper}
{S_1}  \leq CE(S)+C_{d,e}.
\end{align}
To estimate ${S_2}$, we introduce the last multiplier in what follows: \\
\\
Since $(\overline{\Omega \setminus \omega})\cap (\overline{Q_2 \cap \Omega}) = \emptyset$, there exists a function $\beta \in C^\infty_0(\mathbb{R}^2)$ such that 
\begin{equation}\label{beta}
\begin{cases}  
0 \leq \beta \leq 1,
\\ 
\beta = 1\ \ \hbox{on}\ \ Q_2\cap \Omega ,
\\ 
\beta=0\ \ \hbox{on}\ \ \Omega \setminus \omega. 
\end{cases}
\end{equation}
For every $t \geq 0$, let $z$ be the solution of the following elliptic problem:
\begin{equation}\label{ellipt}
\left\lbrace
\begin{array}{ccc}
\Delta z = \beta u & \text{in} \ \Omega , \\
z=0 &\text{on} \ \partial \Omega .
\end{array}
\right.
\end{equation}
One can prove the following lemma:
\begin{Lemma} \label{lemmaellip}
	Under the hypotheses of Proposition \ref{propper} with $z$ as defined in \eqref{ellipt},
	it holds that
	\begin{align}
	&|| z ||_{L^2(\Omega)} \leq C ||u||_{L^2(\Omega)},\ \  
	|| z_t ||^2_{L^2(\Omega)} \leq C\int_\Omega \beta|u_t|^2 \, dx,\ \ 
	\Vert  \nabla z \Vert _{L^2(\Omega)}\leq C ||\nabla u||_{L^2(\Omega)},\label{z1}\\
	\forall &S\leq T \in \mathbb{R_+} ,\ \ \int_S^T \int_\Omega \beta u^2  \, dx\, dt = \left[\int_\Omega zu_t\, dx \right]^T_S + \int_S^T  \int_\Omega \left( -z_tu_t+z\left[ag(u_t+d)+e\right]\right)  \, dx\, dt.\label{id4}
	\end{align}
\end{Lemma} 
\textbf{\emph{Proof of Lemma \ref{lemmaellip}:}}
Equation~\ref{z1} gathers standard elliptic estimates from the definition of $z$ as a solution of \eqref{ellipt} while \eqref{id4} is obtained by using $z$ as a multiplier for $\mathbf{(P_{dis})}$. Steps of the proof are similar to the ones that led to equations $(5.22)$, $(5.25)$ and $(5.26)$ in \cite{Martinez1999}.
\begin{flushright}
	\begin{small}
	$\blacksquare$
\end{small}
\end{flushright}
Since the non negative $\beta$ is equal to $1$ on $Q_2$ and $0$ on $\mathbb{R}^2\setminus \omega$, it follows from \eqref{id4} that
\begin{equation} \label{zineq}
{S_2}  \leq \underbrace{\left[\int_\Omega zu_t\, dx \right]^T_S}_{{U_1}} - \underbrace{\int_S^T \int_\Omega z_tu_t\  \, dx\, dt}_{{U_2}}+\underbrace{\int_S^T \int_\Omega z(ag(u_t+d)+e)\  \, dx\, dt}_{{U_3}}.
\end{equation}
We estimate $U_1$, $U_2$ and $U_3$.
We start by handling ${U_1}$. One has from Cauchy-Schwarz inequality, then \eqref{z1} and Poincaré inequality that
\begin{equation}\label{eq:estU1}
\left|\int_\Omega zu_t\, dx \right| \leq || z ||_{L^2(\Omega)}||u_t ||_{L^2(\Omega)}\leq
C|| \nabla u ||_{L^2(\Omega)}||u_t ||_{L^2(\Omega)}\leq CE(t).
\end{equation}
Using \eqref{eq:estU1} and the fact that $E$ is non-increasing, it is then immediate to derive that
\begin{align}\label{zineq1per}
\left|{U_1} \right|= \left| \left(  \int_\Omega zu_t\, dx \right)(T)- \left(  \int_\Omega zu_t\, dx \right)(S) \right|
\leq C(E(T)+E(S)).
\end{align}
Finally, using \eqref{boundE'} in \eqref{zineq1per}, we obtain that
\begin{align}
{U_1}  \leq CE(S)+C_{d,e}. \label{zu_t}
\end{align}
As for $U_2$, the use of Young inequality with an arbitrary real number $0<\eta_0<1$ yields
$$
\left|{U_2} \right|\leq \int_S^T \int_\Omega \frac{1}{2\eta_0}|z_t|^2  \, dx\, dt + \int_S^T \int_\Omega \frac{\eta_0}{2}|u_t|^2  \, dx\, dt.
$$
Then, we use \eqref{z1} and the fact that $0\leq \beta\leq 1$ to conclude the following estimate:
\begin{align}\label{zhc2per}
{U_2} \leq \frac{C}{\eta_0}\int_S^T  \int_\omega u_t^2  \, dx\, dt+ C\eta_0 \int_S^T E  \, dx\, dt,
\end{align}
where $\eta_0$ is a positive real number to be chosen later.
\\
Left to estimate ${U_3}$. We can rewrite it as the following:
\begin{align}\label{sumV}
U_3 = \underbrace{\int_S^T\int_{|u_t+d|\leq 1} a(x)zg(u_t+d) dx dt}_{V_1} + \underbrace{\int_S^T\int_{|u_t+d|> 1} a(x)zg(u_t+d) dx dt}_{V_2}+ \underbrace{\int_S^T\int_\Omega a(x)ze dx dt}_{V_3}.
\end{align}
We estimate the three terms $V_1$, $V_2$ and $V_3$. We start by estimating $V_1$. We have using Young inequality that
\begin{align}  \label{95}
{V_1} \leq C{\eta_0} \int_S^T E\ dt + \frac{1}{ \eta_0} \int_S^T \int_{|u_t+d|\leq 1} |a g(u_t+d)|^2 \, dx\, dt .
\end{align} 
The fact that $g(0)=0$ implies the existence of a constant $C>0$ such that $|g(x)|\leq C|x|$ for all $|x|\leq 1$. Combining it with the fact that $g(x)x\geq 0, \ \forall\ x \in \mathbb{R} $, it follows that 
\begin{align} \label{ag}
\int_S^T \int_{|u_t+d|\leq 1} |a g(u_t+d)|^2 \, dx\, dt  &\leq   \int_S^T\int_{|u_t+d|\leq 1} a(.)(u_t+d) g(u_t+d) \, dx\,dt  \notag \\
&\leq \int_S^T\int_{\Omega} a(.)(u_t+d) g(u_t+d) \, dx\, dt.
\end{align}
Using \eqref{nonincper} and Young inequality with $0<\eta_1<1$,
\begin{align} 
&\int_S^T\int_{\Omega} a(.)(u_t+d) g(u_t+d) \, dx\,dt = \int_S^T \int_{\Omega} au_t g(u_t+d) \, dx\, dt +\int_S^T\int_{\Omega} ad g(u_t+d) \, dx\, dt \notag
\\
&\leq \int_S^T \int_{\Omega} au_t g(u_t+d) \, dx\, dt + \int_S^T \int_{\Omega} u_t e \, dx\, dt  -\int_S^T \int_{\Omega} u_t e  \, dx\, dt+C\int_S^T\int_{\Omega} |d|| g(u_t+d)| \, dx\, dt \notag \\
& \ \ \ \ \ \ \ \ \ \leq  \int_S^T (-E') dt + \int_S^T \int_{\Omega} |u_t|| e|  \, dx\, dt +C\int_S^T\int_{\Omega} |d|| g(u_t+d)| \, dx\, dt \notag \\ 
&\ \ \ \ \ \ \ \ \ \ \leq E(S)+ C\eta_1 \int_S^T E\ dt +\frac{C}{\eta_1}\int_S^T \int_\Omega| e|^2  \, dx\, dt +C\int_S^T\int_{\Omega} |d|| g(u_t+d)| \, dx\, dt
\notag \\ 
&\ \ \ \ \ \ \ \ \ \ \leq E(S)+ C\eta_1 \int_S^T E\ dt +\frac{C}{\eta_1} \bar{C}_1(e)+C\int_S^T\int_{\Omega} |d|| g(u_t+d)| \, dx\, dt . \label{agg}
\end{align}
Left to estimate $\int_S^T\int_{\Omega} |d|| g(u_t+d)| \, dx\, dt $, we proceed as the following:
\begin{align}
\int_S^T&\int_{\Omega} |d|| g(u_t+d)| \, dx\, dt =	\int_S^T\int_{|u_t+d|\leq 1} |d|| g(u_t+d)| \, dx\, dt +\int_S^T\int_{|u_t+d|> 1} |d|| g(u_t+d)| \, dx\, dt \notag
\\
& \leq C\int_S^T\int_{|u_t+d|\leq 1} |d| \, dx\, dt + \frac{C}{\eta_1'}\int_S^T\int_{|u_t+d|> 1} |d|^2 \, dx\, dt +\eta_1' \int_S^T\int_{|u_t+d|> 1} |g(u_t+d)|^2 \, dx\, dt 
\notag
\\
& \leq CC_6(d) + \frac{C}{\eta_1'} C_1(d) +C\eta_1' \int_S^T\int_{|u_t+d|> 1} |u_t+d|^{2q} \, dx\, dt \notag
\\
& \leq CC_6(d) + \frac{C}{\eta_1'} C_1(d) +C\eta_1' \int_S^T\int_{\Omega} |u_t|^{2q} + C\eta_1'\int_S^T\int_{\Omega} |d|^{2q} \, dx\, dt,
\end{align}
where $0<\eta_1' <1$. Then, using \eqref{gagniper},  
\begin{align}\label{dg}
\int_S^T\int_{\Omega} |d|| g(u_t+d)| \, dx\, dt &\leq CC_6(d) + \frac{C}{\eta_1'} C_1(d) +\eta_1'(C_u+C_{d,e}) \int_S^T
E(t)\, dt+ C\eta'_1 C_1(d)\notag
\\
&\leq  \frac{1}{\eta_1'} C_{d,e} +\eta_1'(C_u+C_{d,e}) \int_S^T
E(t)\, dt. 
\end{align}
Combining \eqref{agg} and \eqref{dg}, 
\begin{align}\label{ugu}
\int_S^T\int_{\Omega} a(.)(u_t+d) g(u_t+d) \, dx\,dt \leq E(S)+\left(\eta_1 +\eta_1'(C_u+C_{d,e})\right)\int_S^T E\ dt + \frac{1}{\eta_1 \eta_1'}C_{d,e}.
\end{align}
Combining now \eqref{ugu}, \eqref{agg} and \eqref{95}, we obtain that
\begin{align*}
{V_1}
  &\leq C\left( \eta_0 + \frac{\eta_1'}{\eta_0}(C_u+C_{d,e})+\frac{\eta_1}{\eta_0}\right)\int_S^T E\ dt + \frac{C}{\eta_0} E(S)+ \frac{1}{\eta_1\eta_0 \eta_1'}C_{d,e}.
\end{align*}
We take $\eta_1= \eta_0^2$ and $\eta_1'= \frac{\eta_0^2}{C_u+C_{d,e}}$ if 
$C_u+C_{d,e}>0$. In that case, $V_1$ would be estimated by
\begin{align} \label{zhg1per}
{V_1}
&\leq C\eta_0\int_S^T E\ dt + \frac{C}{\eta_0} E(S)+ \frac{1}{\eta_0^5}C_{d,e}(C_u+C_{d,e}).
\end{align}   
If $C_u=C_{d,e}=0$, the above equation holds true trivially.
\begin{Remark}
With such a choice of $\eta_1$ and $\eta_1'$, we have the following useful estimate obtained from \eqref{ugu}:
\begin{align}\label{useful}
\int_S^T\int_{\Omega} a(.)(u_t+d) g(u_t+d) \, dx\,dt \leq E(S)+C\eta_0^2\int_S^T E\ dt + \frac{1}{\eta_0^4}(C_{d,e}C_u+C_{d,e}).
\end{align}
\end{Remark}
To estimate ${V_2}$, first notice that from Rellich-Kondrachov's theorem in dimension 
two (cf. \cite{Brezis}) that
$
H^1(\Omega) \subset L^{q+1}(\Omega),
$
which means that $\exists \ C>0 \ $ such that
$
\Vert  z \Vert _{L^{q+1}(\Omega)} \leq C \Vert  z \Vert _{H^1(\Omega)}
$, adding to that the fact that $z\in H^1_0(\Omega) $ and  \eqref{z1}, it holds that
\begin{align}\label{embres}
\Vert  z \Vert _{L^{q+1}(\Omega)} \leq C \sqrt{E}.
\end{align}
Then, using H\"older inequality yields
	\begin{align}\label{h1}
{V_2} \leq& \int_S^T \left(\int_{|u_t+d|>1}  (a|g(u_t+d)|)^{\frac{q+1}{q}} \, dx \right) ^{\frac{q}{q+1}} \left(\int_{|u_t+d|>1}  |z|^{q+1} \, dx \right) ^{\frac{1}{q+1}} dt.
\end{align}
Combining \eqref{h1} with the hypothesis given by \eqref{hypgq}, we get that
	\begin{align*}
{V_2} &\leq C \int_S^T \left(\int_{|u_t+d|>1}  a|u_t+d|| g(u_t+d)| \, dx \right) ^{\frac{q}{q+1}} \left(\int_{|u_t+d|>1}  |z|^{q+1} \, dx \right) ^{\frac{1}{q+1}} dt.
\end{align*}
Using Young inequality for an arbitrary $0<\eta_2<1$,
	\begin{align}\label{h2}
V_2& \leq C \int_S^T \left(\frac{1}{\eta_2^{\frac{q+1}{q}}}\int_{|u_t+d|>1}  a(x)(u_t+d)g(u_t+d) \, dx  +\eta_2^{q+1} \int_{\Omega}  |z|^{q+1} \, dx \right) dt
\notag
\\
&\leq C \int_S^T \left(\frac{1}{\eta_2^{\frac{q+1}{q}}}\int_{\Omega}  a(x)(u_t+d)g(u_t+d)\, dx  +\eta_2^{q+1} \int_{\Omega}  |z|^{q+1} \, dx \right) dt
\notag
\\
&\leq C \int_S^T \left(\frac{1}{\eta_2^{\frac{q+1}{q}}}\int_{\Omega}  a(x)u_tg(u_t+d)\, dx  + \frac{C}{\eta_2^{\frac{q+1}{q}}}\int_{\Omega}  |d||g(u_t+d)|\, dx  + \eta_2^{q+1} \int_{\Omega}  |z|^{q+1} \, dx \right) dt.
\end{align}
The previous inequality combined with \eqref{nonincper} and \eqref{embres} implies that
\begin{align*}
V_2
&\leq C \int_S^T \left(\frac{1}{\eta_2^{\frac{q+1}{q}}} (-E') -\frac{1}{\eta_2^{\frac{q+1}{q}}}\int_\Omega u_te \, dx +\frac{C}{\eta_2^{\frac{q+1}{q}}}\int_{\Omega}  |d||g(u_t+d)|\, dx  +\eta_2^{q+1} E^{\frac{q+1}{2}}\right) dt.
\end{align*}
Then, using \eqref{boundE'}, $E$ satisfies
\begin{align}
	\int_S^T E^{\frac{q+1}{2}} \,dt = \int_S^T E^{\frac{q-1}{2}}E \, dt \leq  (CE(0)+C_{d,e})^{\frac{q-1}{2}} \int_S^T E \, dt  \leq  (C_u+C_{d,e}) \int_S^T E \, dt,
\end{align}
which gives that
\begin{align*}
V_2 &\leq \frac{C}{\eta_2^{\frac{q+1}{q}}} E(S) +\eta_2 ^{q+1} (C_u+C_{d,e}) \int_S^T E  dt 
 + \frac{C}{\eta_2^{\frac{q+1}{q}}} \int_S^T \left(-\int_\Omega u_te \, dx + \int_{\Omega}  |d||g(u_t+d)|\, dx \right) dt.
\end{align*}
We fix $\eta_2 = \left(\frac{\eta_0}{(C_u+C_{d,e})} \right)^{\frac{1}{q+1}}$. It follows that 
\begin{align*}
	\eta_2^{q+1} (C_u+C_{d,e})&= \eta_0, \\
	\frac{C}{\eta_2^{\frac{q+1}{q}}}= C \frac{(C_u + C_{d,e})^\frac{1}{q}}{\eta_0 ^\frac{1}{q}} \leq \frac{C}{\eta_0}& (C_u^\frac{1}{q} + C_{d,e}^\frac{1}{q}) = \frac{1}{\eta_0^\frac{1}{q}}\left(C_u + C_{d,e}\right),
\end{align*}
which leads to 
\begin{small}
	\begin{align}
{V_2} \leq \frac{1}{\eta_0^\frac{1}{q}}(C_u + C_{d,e}) E(S) +\eta_0 \int_S^T E  dt 
+ \frac{1}{\eta_0^\frac{1}{q}}(C_u + C_{d,e})\int_S^T \left(-\int_\Omega u_te \, dx + \int_{\Omega}  |d||g(u_t+d)|\, dx \right) dt. \label{zppp1}
\end{align}
\end{small}
\ \\
To finish the estimation of ${V_2}$, we still have to handle the last two integral terms in \eqref{zppp1}.
\\
On one hand, we have already estimated the term $ \int_S^T \int_{\Omega}|d||g(u_t+d)|\,dx\, dt$ in \eqref{dg}. We have immediately for some $0<\eta_3<1$ that
\begin{small}
	\begin{align}
(C_u+C_{d,e})\int_S^T \int_\Omega  |d||g(u_t+d)|\, dx\, dt 
&\leq  \eta_3 (C_u+C_{d,e}) \int_S^T E\ dt  + \frac{1}{\eta_3}(C_{d,e}C_u+C_{d,e}).
\end{align}
\end{small}
Choosing $\eta_3$ to be equal to $\frac{\eta_0^{\frac{q+1}{q}}}{(C_u+C_{d,e})}$ implies that
\begin{align*}
	\eta_3(C_u+C_{d,e}) &= \eta_0^\frac{q+1}{q},\\ \frac{1}{\eta_3}  (C_{d,e}C_u+C_{d,e})&\leq \frac{1}{\eta_0^\frac{q+1}{q}}(C_{d,e}C_u+C_{d,e}) ,
\end{align*}
which gives that
\begin{align}\label{zppp2}
(C_u+C_{d,e})\int_S^T \int_\Omega |d||g(u_t+d)|\, dx\, dt
\leq \eta_0^\frac{q+1}{q} \int_S^T E\ dt  + \frac{1}{\eta_0^\frac{q+1}{q}}\left(C_{d,e}C_u+ C_{d,e}\right).
\end{align}
On the other hand, we have for $0<\eta_4<1$ that
\begin{align*}
(C_u+C_{d,e})\int_S^T \int_\Omega u_te  \, dx\, dt \leq \eta_4 (C_u+C_{d,e})\int_S^T E\ dt + \frac{1}{\eta_4} (C_{d,e}C_u+ C_{d,e}).
\end{align*}
Using the same concept as before, we fix $\eta_4=\frac{\eta_0^\frac{q+1}{q}}{C_u+C_{d,e}}$, we obtain that
\begin{align}\label{u_td2est}
(C_u+C_{d,e})\int_S^T \int_\Omega u_te  \, dx\, dt \leq \eta_0^\frac{q+1}{q} \int_S^T E\ dt + \frac{1}{\eta_0^\frac{q+1}{q}}\left(C_{d,e}C_u+ C_{d,e}\right),
\end{align}
Combining \eqref{zppp1}, \eqref{zppp2} and \eqref{u_td2est}, we conclude that the estimation of $V_2$ is given by
\begin{align}\label{zppp3}
{V_2}\leq \frac{1}{\eta_0^\frac{1}{q}}(C_u+C_{d,e}) E(S) &+\eta_0\int_S^T E  dt +\frac{1}{\eta_0^\frac{q+2}{q}}\left( C_{d,e}C_u+ C_{d,e}\right).
\end{align}
As for $V_3$, we simply have when using \eqref{z1} and Young inequality with $\eta_0$ that
\begin{align*}
{V_3} \leq C\eta_0 \int_S^T E\ dt +\frac{C}{\eta_0} \bar{C}_1(e),
\end{align*}
which means that
\begin{align}\label{zd2}
	V_3 \leq C\eta_0 \int_S^T E\ dt + \frac{1}{\eta_0}  C_{d,e} .
\end{align}
To achieve an estimation of ${S_2}$, we just combine  \eqref{zu_t},\eqref{zhc2per} \eqref{zhg1per},  \eqref{zppp3} and \eqref{zd2} to get 
\begin{align}\label{prev}
{S_2}  &\leq  C\eta_0 \int_S^T E\ dt +  \frac{C}{\eta_0} \int_S^T\int_\omega u_t^2  \, dx\, dt +\left(C+\frac{C}{\eta_0}+\frac{1}{\eta_0^\frac{1}{q}}(C_u+C_{d,e})\right) E(S) 
\notag
\\
&\ \ \ \ \ \ \ \ + \left(\frac{1}{\eta_0^5}+\frac{1}{\eta_0^\frac{q+2}{q}}\right)\left(C_{d,e}C_u+ C_{d,e}\right)+ \left(\frac{1}{\eta_0}+1\right)C_{d,e}.
\end{align} 
We can simplify the previous estimate by using the fact that $ 0<\eta_0<1$. As a result, \eqref{prev} becomes
\begin{align}\label{hc1per}
{S_2}  &\leq  C\eta_0 \int_S^T E\ dt +  \frac{C}{\eta_0} \int_S^T\int_\omega u_t^2  \, dx\, dt +\frac{1}{\eta_0}\left(C+C_u+C_{d,e}\right) E(S) 
+ \frac{1}{\eta_0^5}\left(C_{d,e}C_u+ C_{d,e}\right).
\end{align} 
Regarding ${S_3}$, we follow the same steps we followed to get ${V_1}+{V_2}$. It is possible because $u$ satisfies the same result \eqref{embres} as $z$ from before. Hence, we obtain that
\begin{align}\label{hcper} 
{S_3}&\leq   C\eta_0 \int_S^T E\ dt  +\frac{1}{\eta_0}\left(C+C_u+C_{d,e}\right) E(S) + \frac{1}{\eta_0^5}\left(C_{d,e}C_u+ C_{d,e}\right).
\end{align}
Finally, to estimate $S_4$, we simply have when using young inequality that
\begin{align}\label{ud2}
{S_4} \leq \eta_0 \int_S^T E\ dt + \frac{1}{\eta_0}C_{d,e}.
\end{align}
We complete the estimate of ${T_2}$ in \eqref{Aper} by combining the estimations of $S_1$, $S_2$, $S_3$ and $S_4$.
Hence the proof of Lemma \ref{es2per} is completed.
\begin{flushright}
\small{$\blacksquare$}
\end{flushright}
An estimate of ${T_3}$ is provided in the next lemma:
\begin{Lemma}\label{es3per}
Under the hypotheses of Proposition \ref{propper}, we have the following estimate:
\begin{small}
	\begin{align} \label{estim2per}
{T_3} \leq C\eta_0 \int_S^T E\ dt  +\frac{1}{\eta_0}\left[C+(1+C_{\eta_{0}})(C_u+C_{d,e})\right] E(S) + \frac{1}{\eta_0^5}\left(C_{\eta_{0}}^3+1\right)\left(C_{d,e}C_u+ C_{d,e}\right),
\end{align}
\end{small}
where $0<\eta_0<1$  is a positive arbitrary real number to be chosen later and $C_{\eta_{0}}$ is an implicit positive constant that depends on $\eta_{0}$ only.
\end{Lemma}
\textbf{\emph{Proof of Lemma~\ref{es3per}:}} First, note that 
\begin{equation}\label{Mu}
{T_3} \leq \frac{1}{2}{S_3} + \underbrace{\int_S^T \int_{\Omega} |ag(u_t+d) \nabla u. h| \, dx\, dt}_{{X}}.
\end{equation}
We have already estimated ${S_3}$ in \eqref{hcper}. It remains to deal with ${X}$. 
Using Young inequality implies that
\begin{align}
{X} &\leq \frac{C}{\eta_0}\int_S^T \int_{\Omega}  (a|g(u_t+d)|)^2 \, dx\, dt+C\eta_0 \int_S^T \int_{\Omega}  |\nabla u|^{2} \, dx\, dt
\notag
\\
&\leq \frac{C}{\eta_0}\int_S^T \int_{\Omega}  a|g(u_t+d)|^2 \, dx\, dt+C\eta_0 \int_S^T E\ dt. \label{nug111per}
\end{align}
Now, set $R_1>1$ to be chosen later. We can rewrite the term $\int_S^T \int_{\Omega}  a|g(u_t+d)|^2 \, dx\, dt$ as
\begin{align}\label{y1y2}
\int_S^T \int_{\Omega}  a|g(u_t+d)|^2  \, dx\, dt=  \underbrace{\int_S^T \int_{|u_t+d|\leq R_1}  a|g(u_t+d)|^2  \, dx\, dt}_{{Y_1}} + 	 \underbrace{\int_S^T \int_{|u_t+d|> R_1}  a|g(u_t+d)|^2  \, dx\, dt}_{{Y_2}}.
\end{align}
Since $g(0)=0$, it holds that $|g(x)| \leq C_{R_1}|x|$ for some constant $C_{R_1}$ and for $|x|<R_1$. Combine it with \eqref{useful}, it follows that $Y_1$ satisfies for some $0<\eta_5<1$
\begin{align} 
{Y_1}
&\leq C_{R_1} \int_S^T \int_{|u_t+d|\leq R_1} |a g(u_t+d)||u_t+d| \, dx\, dt \notag
\\
&\leq C_{R_1} \int_S^T \int_{\Omega} |a g(u_t+d)||u_t+d| \, dx\, dt
\notag
\\
&\leq  C_{R_1}E(S)+C C_{R_1}\eta_5^2\int_S^T E\ dt + \frac{ C_{R_1}}{\eta_5^4}(C_{d,e}C_u+C_{d,e}).
\end{align}
Taking $\eta_5=\frac{\eta_{0}}{\sqrt{C C_{R_1}}}$ leads to
\begin{align}\label{leqR}
	Y_1 \leq  C_{R_1}E(S)+\eta_0^2\int_S^T E\ dt + \frac{ C_{R_1}^3}{\eta_0^4}(C_{d,e}C_u+C_{d,e})
\end{align} 
As for $Y_2$, we use \eqref{hypgq} to obtain that
\begin{align*}
	{Y_2}&\leq C 	\int_S^T \int_{|u_t+d|> R_1} |u_t+d|^{2q} \, dx\, dt \notag
	\\
	&\leq C 	\int_S^T \int_{|u_t+d|> R_1} |u_t|^{2q}  \, dx\, dt +C	\int_S^T \int_{|u_t+d|> R_1} |d|^{2q} \, dx  \notag
	\\
	&\leq C 	\int_S^T \int_{|u_t+d|> R_1} \frac{|u_t+d|}{R_1}|u_t|^{2q}  \, dx\, dt +C	\int_S^T \int_{\Omega} |d|^{2q} \, dx\, dt
	\notag
	\\&\leq C 	\int_S^T \int_{\Omega} \frac{|u_t|}{R_1}|u_t|^{2q}  \, dx\, dt +  C 	\int_S^T \int_{\Omega} \frac{|d|}{R_1}|u_t|^{2q}  \, dx\, dt+	CC_1(d) \notag
	\\
	&\leq \frac{C}{R_1} 	\int_S^T \int_{\Omega} |u_t|^{2q+1}  \, dx\, dt +\frac{C}{R_1^2}\int_S^T \int_{\Omega}|u_t|^{4q}  \, dx\, dt+C_{d,e}.	
\end{align*}
Then, we use Lemma \ref{gagnircorper} as well as the fact that $R_1>1$ to conclude that $Y_2$ satisfies
\begin{align*}
{Y_2}
&\leq \frac{1}{R_1} (C_u + C_{d,e})	\int_S^T E\ dt +	C_{d,e}.
\end{align*}
We take $R_1 = \frac{(C_u+C_{d,e})}{\eta_{0}^2}$, we get the simplified estimate
\begin{align}
{Y_2}
&\leq \eta_{0}^2	\int_S^T E\ dt +	C_{d,e}.	\label{geqR}
\end{align}
\begin{Remark}
	For such a choice of $R_1$, and based on how $C_{R_1}$ is defined, we can assume that $C_{R_1}$ in \eqref{leqR} is a constant of the type $C_{\eta_0}(C_u+C_{d,e})$, where $C_{\eta_0}$ is a positive constant that depends on $\eta_0$ only.
\end{Remark}
Combining \eqref{nug111per}, \eqref{y1y2}, \eqref{leqR} and \eqref{geqR} implies that
\begin{small}
	\begin{align}\label{nhcper2}
{X} &\leq C\eta_0 \int_S^T E\ dt  +{\frac{C_{\eta_0}}{\eta_{0}}}(C_{d,e}+C_u) E(S)+ \frac{C^3_{\eta_0}}{\eta_{0}^5}(C_{d,e}C_u+ C_{d,e})+ \frac{C_{d,e}}{\eta_{0}}.
\end{align}
\end{small}
Finally, we combine \eqref{Mu} and \eqref{nhcper2} with the estimation of  ${S_3}$, we obtain \eqref{estim2per}.
\begin{flushright}
\small{$\blacksquare$}
\end{flushright}
We next seek to prove the upper bound of ${T_4}$ that is given by the following lemma
\begin{Lemma}\label{lemmaT4}
	Under the hypotheses of Proposition \ref{propper}, the following estimate holds:
\begin{align}\label{estd2m}
{T_4} \leq C\eta_0 \int_S^T E\ dt +  \frac{C}{\eta_0} C_{d,e},
\end{align}
where $0<\eta_{0}<1$ is a positive constant to be chosen later.
\end{Lemma}
\textbf{\emph{Proof of Lemma \ref{lemmaT4}}:} We have that
\begin{equation}\label{hand}
{T_4}\leq \frac{1}{2}\int_S^T  \int_{\Omega} |eu| \, dx\, dt + \int_S^T \int_{\Omega} |e \nabla u. h |\, dx\, dt.
\end{equation}
On one hand, using Young inequality gives that
\begin{align}\label{hand1}
	\int_S^T  \int_{\Omega} |eu| \, dx\, dt \leq \eta_{0} \int_S^T E\ dt + \frac{C}{\eta_0} C_{d,e}.
\end{align}
On the other hand, it gives that
\begin{align}\label{hand2}
	\int_S^T \int_{\Omega} |e \nabla u. h |\, dx\, dt \leq  \eta_{0} \int_S^T E\ dt +  \frac{C}{\eta_0} C_{d,e}
\end{align}
Combining \eqref{hand}, \eqref{hand1} and \eqref{hand2}, we prove \eqref{estd2m}.
\begin{flushright}
	\small{$\blacksquare$}
\end{flushright}
It remains to handle the last term ${T_5}$. 
\begin{Lemma}\label{lemmaT5}
	Under the hypotheses of Proposition \ref{propper}, we have the following estimation:
	\begin{align}\label{estu_t}
	{T_5} \leq \eta_{0} \int_S^T E\ dt +\bar{C}_{\eta_{0}} (C_u+C_{d,e})E(S)+\frac{\bar{C}_{\eta_{0}}^3}{\eta_0^2} \left(C_{d,e}C_u+C_{d,e}\right)+C_{d,e},
	\end{align}
	where $0<\eta_{0}<1$ is a positive constant to be chosen later and and $\bar{C}_{\eta_{0}}$ is an implicit positive constant that depends on $\eta_{0}$ only.
\end{Lemma}
\textbf{\emph{Proof of Lemma \ref{lemmaT5}}:}
For every $R_2> 1$, we have that
\begin{align}
{T_5} &\leq \frac{1}{a_0} \int_S^T \int_\omega a(x) u_t^2 \, dx\, dt \leq C \int_S^T \int_\Omega a(x)(u_t+d)^2 \, dx\, dt + C \int_S^T \int_\Omega a(x) d^2 \, dx\, dt \notag
\\
 &\leq C \underbrace{\int_S^T \int_{|u_t+d|\leq R_2} a(x) (u_t+d)^2 \, dx\, dt}_{{Z_1}} +C\underbrace{\int_S^T \int_{|u_t+d|>R_2} a(x) (u_t+d)^2 \, dx\, dt}_{{Z_2}} + CC_1(d). \label{*0}  
\end{align}
On one hand, since $g'(0)>0$, there exists $\alpha_{R_2}>0$ such that 
$|g(v)| \geq \alpha_{R_2} |v|$ \ for $|v| \leq R_2$. Combining that with \eqref{useful} yields for some $0<\eta_6<1$
\begin{align*}
 {Z_1} &\leq  \int_S^T \int_{|u_t+d|\leq R_2} a(x) (u_t+d) g(u_t+d)\frac{(u_t+d)}{g(u_t+d)} \, dx\, dt \notag
  \\
 & \leq \frac{1}{\alpha_{R_2}} \int_S^T \int_{|u_t+d|\leq R_2} a(x) (u_t+d) g(u_t+d) \, dx\, dt   \notag 
 \\
  & \leq  \frac{1}{\alpha_{R_2}} \int_S^T \int_{\Omega} a(x) (u_t+d) g(u_t+d) \, dx\, dt  \notag 
  \\
  & \leq \frac{1}{\alpha_{R_2}}E(S)+C\frac{1}{\alpha_{R_2}}\eta_6^2\int_S^T E\ dt +\frac{1}{\alpha_{R_2}} \frac{1}{\eta_6^4}(C_{d,e}C_u+C_{d,e}).
\end{align*}
We choose $\eta_6=\sqrt{\frac{{\alpha_{R_2}}}{{C }}\eta_{0}}$, we obtain that
\begin{align*}
{Z_1} \leq \frac{1}{\alpha_{R_2}} E(S)+\eta_{0} \int_S^T E\ dt +\frac{1}{\alpha_{R_2}^3\eta_{0}^2}(C_{d,e}C_u+C_{d,e}). 
\end{align*}
As for $Z_2$, we have that 
\begin{align}{Z_2}
&\leq C \int_S^T \int_{|u_t+d|>R_2} |u_t|^{2} \, dx\, dt + C\int_S^T \int_{|u_t+d|>R_2} |d|^{2}  \, dx\, dt  \notag 
\\
&\leq C \int_S^T \int_{|u_t+d|>R_2} \frac{|u_t+d|}{R_2}|u_t|^{2} \, dx\, dt +CC_1(d) \notag
\\
&\leq C \int_S^T \int_{|u_t+d|>R_2} \frac{|u_t|^3}{R_2} \, dx\, dt + C \int_S^T \int_{|u_t+d|>R_2} \frac{|u_t|^2|d|}{R_2} \, dx\, dt +CC_1(d) \notag \\
&\leq \frac{C}{R_2} \int_S^T \int_{|u_t+d|>R_2} |u_t|^3\, dx\, dt + \frac{C}{R_2^2}\int_S^T \int_{|u_t+d|>R_2}  |u_t|^4 \, dx\, dt +CC_1(d).
\label{Rin}
\end{align}
We use Lemma \ref{gagnircorper} and the fact that $R_2>1$, we derive the following:
\begin{align}
\frac{C}{R_2}& \int_S^T \int_{|u_t+d|>R_2} |u_t|^3\, dx\, dt + \frac{C}{R_2^2}\int_S^T 
\int_{|u_t+d|>R_2}  |u_t|^4 \, dx\, dt \leq  \left(\frac{C_u+C_{d,e}}{R_2}    \right) \int_S^T E\ dt. \label{Rin1}
\end{align}
We choose $R_2= \frac{(C_u+C_{d,e})}{\eta_{0}}$ and we 
combine \eqref{Rin} and \eqref{Rin1} we have that
\begin{align}
{Z_2} \leq  \eta_{0} \int_S^T E\ dt +C_{d,e}. \label{2*}
\end{align}
\begin{Remark}
For such a choice of $R_2$, and based on how $\alpha_{R_2}$ is defined, we can 
assume that $\frac{1}{\alpha_{R_2}}$ is also a constant of the type $\bar{C}_{\eta_{0}} \left(C_u+C_{d,e}\right)
$, where $\bar{C}_{\eta_{0}}$ is a constant that depends on $\eta_{0}$ only. As a result, $Z_1$ is estimated by
\begin{align}
{Z_1} 
\leq \eta_{0} \int_S^T E\ dt +\bar{C}_{\eta_{0}} (C_u+C_{d,e})E(S)+\frac{\bar{C}_{\eta_{0}}^3}{\eta_0^2} \left(C_{d,e}C_u+C_{d,e}\right). \label{*1}
\end{align}
\end{Remark}
Combining \eqref{*0}, \eqref{2*} and \eqref{*1} and using \eqref{bornc1} and \eqref{lip2}, it follows that  
\begin{align*}
{T_5} \leq \eta_{0} \int_S^T E\ dt +\bar{C}_{\eta_{0}} (C_u+C_{d,e})E(S)+\frac{\bar{C}_{\eta_{0}}^3}{\eta_0^2} \left(C_{d,e}C_u+C_{d,e}\right)+C_{d,e},
\end{align*}
which proves Lemma \ref{lemmaT5}.
\begin{flushright}
	\small{$\blacksquare$}
\end{flushright}
The estimation of $T_5$ gives a direct estimation of the term $\frac{C}{\eta_0}\int_S^T  \int_\omega u_t^2  \, dx\, dt$ left in the estimation of $T_2$. We can easily manage to have that
\begin{align}
\frac{1}{\eta_0}\int_S^T  \int_\omega u_t^2  \, dx\, dt \leq \eta_{0} \int_S^T E\ dt +\frac{1}{\eta_0}\bar{C}_{\eta_{0}^2} (C_u+C_{d,e})E(S)+\frac{1}{\eta_0}\frac{\bar{C}_{\eta_{0}^2}^3}{\eta_0^4} \left(C_{d,e}C_u+C_{d,e}\right)+C_{d,e}.
\end{align}
It is obtained by following the same steps that led to the estimation of $T_5$ with replacing $\eta_0$ by $\eta_0^2$.
\\ \\
We can finally finish the proof of Proposition \ref{propper}: we combine the estimations of $T_i, \ i=1,2,3,4,5$, which are given by
\eqref{estim1per},  \eqref{ineqnablaper},  \eqref{estim2per}, \eqref{estd2m} and 
\eqref{estu_t} with \eqref{ineq2per}, then
we choose $\eta_0$ such that $C\eta_{0}<1$, which means that the term $C\eta_{0} \int_S^T E(t)\, dt$  gets absorbed by $\int_S^T E(t)\, dt$. Then we use the fact that $C_{d,e}E(S)\leq C_{d,e}(E(0) + C_{d,e})=C_{d,e}C_u + C_{d,e} $ and the fact that
the choice of $\eta_{0}$ will be a constant $C$, we obtain \eqref{lem73per}.
\begin{flushright}
	\small{$\blacksquare$}
\end{flushright}
\textbf{\emph{Proof of the energy estimate of Theorem~\ref{thrmdis}:}}
Using the key result given by \eqref{lem73per}, we get at once from 
Theorem~\ref{gronw1lem} that \eqref{gronw0} holds true with $T=C+C_u$ and 
$C_0=(1+C_u)C_{d,e}$. Using \eqref{boundE'} for $t\geq 1$ with $T=t$ and 
$S\in [t-1,t]$ and integrating it over $[t-1,t]$, one gets that 
$$
E(t)\leq C\int_{t-1}^tE(s)\ ds+C_{d,e}\leq C\int_{t-1}^{\infty}E(s)\ ds+C_{d,e}.
$$
Combining the above with \eqref{hyp} yields \eqref{gt1} for $t\geq 1$. In turn, 
\eqref{boundE'} with $T\in [0,1]$ and $S=0$ provides \eqref{gt1} for $t\leq 1$. The proof of Theorem~\ref{thrmdis} is then completed. 
\begin{flushright}
	\small{$\blacksquare$}
\end{flushright}
\appendix 
\addcontentsline{toc}{section}{Appendix}
\section{Appendix}\label{background}
We list in what follows, technical results used in the core of the paper.
\begin{Theorem}{\textbf{Gronwall integral lemma }} \label{gronw1lem} \\
Let $E:\mathbb{R_+}\to\mathbb{R}_+$ satisfy for some $C_0,T>0$:
\begin{equation}\label{hyp}
\int_t^{+\infty} E(s)ds \leq TE(t)+ C_0,\ \ \forall\ t \geq 0.
\end{equation}
Then, the following estimate hold true
\begin{equation}\label{gronw0}
\int_t^{+\infty} E(s)ds \leq TE(0) e^{-\frac{t}{T}}   + C_0,\ \ \forall\ t \geq 0.
\end{equation}
If in addition, $t\mapsto E(t)$ is non-increasing, one has
\begin{equation}\label{gronw1}
E(t) \leq E(0) e^{1-\frac{t}{T}}   + \frac{C_0}{T}, \ \ \forall\ t \geq 0. 
\end{equation}
\end{Theorem}
The proof is classical, cf. for instance \cite{Alabau}.

\begin{Theorem}{\textbf{Generalized Gronwall lemma}} \label{gronw2lem} \\
Let $F,h_1$ and $h_2$ non negative functions defined on $\mathbb{R}_+$
satisfying
$$
\Vert h_1\Vert_1:=\int_0^\infty h_1(t)dt <\infty,\quad  \Vert h_2\Vert_1:=\int_0^\infty h_2(t)dt <\infty,
$$ 
and
\begin{equation}\label{hyp2}
F(T) \leq F(S)+C_3+ C_1\int_S^{T} h_1(s)F^{\alpha_1}(s)ds + C_2 \int_S^{T}
h_2(s)F^{\alpha_2}(s)ds,\ \ \forall\ S\leq T,
\end{equation}
where $C_1,C_2,C_3$ are positive constants and $0\leq \alpha_1,\alpha_2 <1$. 
Then, $F$ satisfies the following bound
\begin{equation}\label{gronw2}
\sup_{t\in [S,T]}F(t) \leq \max\Big(2(F(S)+C_3),(2\tilde{C})^{\frac{1}{1-\alpha}}\Big),
\hbox{ with }\ \tilde{C}:=C_1\Vert h_1\Vert_1+C_2\Vert h_2\Vert_1,
\end{equation}
where $\alpha:=\max(\alpha_1,\alpha_2)$ if $2\tilde{C}\geq 1$ or 
$\alpha:=\min(\alpha_1,\alpha_2)$ if $2\tilde{C}< 1$.
\end{Theorem}
\textbf{\emph{Proof of Theorem~\ref{gronw2lem}:}} 
Fix $T\geq S\geq 0$. For $t\in [S,T]$ set 
$Y(t)$ for the right-hand side of \eqref{hyp2} applied at the pair of times $S\leq t$. It 
defines a non decreasing absolutely continuous function. Since $F(t)\leq Y(t)\leq 
Y(T)$ for $t\in [S,T]$, one deduces that $F_{S,T}:=\sup_{t\in[S,T]}F(t)$ is finite for every 
$t\in [S,T]$. One gets from \eqref{hyp2} that
$$
F_{S,T}\leq F(S)+ C_3+\tilde{C}\max(F_{S,T}^{\alpha_1},F_{S,T}^{\alpha_2}),
$$
with the notations of \eqref{gronw2}. The latter follows at once by considering 
whether $F(S)+C_3>\tilde{C}\max(F_{S,T}^{\alpha_1},F_{S,T}^{\alpha_2})$ or not.
\\
\\
We recall the following useful result, cf. for instance \cite{Martinez2000}.
\begin{Theorem}{\textbf{Gagliardo–Nirenberg interpolation inequality}}\label{gagnir} \\
Let $\Omega\subset\mathbb{R}^N$ be a bounded Lipschitz domain,  $N\geq 1$, $1\leq r<p \leq \infty, 1\leq q \leq p$ and $m\geq 0$. Then the inequality
\begin{equation}
\Vert  v \Vert _p \leq C \Vert  v \Vert _{m,q}^\theta \Vert  v \Vert _r^{1-\theta} \ \ \ \ for \ \ \ v\in W^{m,q}(\Omega)\cap L^r(\Omega)
\end{equation}
holds for some constant $C>0$ and 
\begin{equation}
\theta = \left(\frac{1}{r}-\frac{1}{p} \right) \left(\frac{m}{N}+\frac{1}{r}-\frac{1}{q} \right)^{-1},
\end{equation}
where $0<\theta\leq 1\ \ (0<\theta < 1$ if $p=\infty$ and $mq=N)$ and $\Vert  . \Vert _p$ denotes the usual $L^p(\Omega)$ norm and $\Vert  . \Vert _{m,q} $ the norm in $W^{m,q}(\Omega)$.
\end{Theorem}
\bibliographystyle{plain} 
\bibliography{biblio} 	
	
\end{document}